\input amstex
\input amsppt.sty
\def\phi{\varphi}
\def\epsilon{\varepsilon}
\NoBlackBoxes

\topmatter
\title
Asymptotic behaviour of the sectional curvature of the Bergman metric for annuli
\endtitle

\author
W\l odzimierz Zwonek
\endauthor

\abstract
We extend and simplify results of \cite{Din~2009} where the asymptotic behavior of the holomorphic sectional curvature of the Bergman metric in annuli is studied. Similarly as in \cite{Din~2009} the description enables us to construct an infinitely connected planar domain (in our paper it is a Zalcman type domain) for which the supremum of the holomorphic sectional curvature is two whereas its infimum is equal to $-\infty$.
\endabstract

\address
Instytut Matematyki, Uniwersytet Jagiello\'nski, \L ojasiewicza 6, 30-348 Krak\'ow, Poland
\endaddress
\email
Wlodzimierz.Zwonek\@im.uj.edu.pl
\endemail
\thanks The research was partially supported by the Research Grant No. N N201 361436 of the Polish Ministry of Science
and Higher Education.
\endthanks
\thanks
2000 Mathematics Subject Classification.
Primary: 32F45. Secondary: 32A36, 30C40
\endthanks
\endtopmatter
\document
For a domain $D\subset\Bbb C^n$, $j=0,1,\ldots$, $z\in D$, $X\in\Bbb C^n$ define
$$
\multline J_D^{(j)}(z;X):=\\
\sup\{|f^{(j)}(z)(X)|^2:f\in
L_h^2(D),f(z)=0,\ldots,f^{(j-1)}(z)=0, ||f||_{L^2(D)}\leq 1\}.
\endmultline
$$
Note that the functions above are the squares of operator norms of
continuous operators defined on a closed subspace of $L_h^2(D)$. 

Let us restrict ourselves to the case when $D$ is bounded. Note that
$J_D^{(0)}(z;X)$ is independent of $X\neq 0$ and is equal to the
Bergman kernel $K_D(z,z)$. Moreover, we may express the Bergman
metric as $\beta_D^2(z;X)=\frac{J_D^{(1)}(z;X)}{J_D^{(0)}(z;X)}$,
$X\neq 0$. And finally the sectional curvature is given by the
formula
$$
R_D(z;X)=2-\frac{J_D^{(0)}(z;X)J_D^{(2)}(z;X)}{J_D^{(1)}(z;X)^2},\;X\neq 0.
$$
Below we list a number of simple properties of the above functions.



The transformation formula for a biholomorphic mapping $F:D_1\mapsto D_2$ is the following
$$
J_{D_1}^{(j)}(z;X)=|\det
F^{\prime}(z)|^{2}J_{D_2}^{(j)}(F(z);F^{\prime}(z)X), $$ from
which we get, among others, the independence of the sectional curvature for
biholomorphic mappings $R_{D_1}(z;X)=R_{D_2}(F(z);F^{\prime}(z)X)$.

If $D_1\subset D_2$ then $J_{D_1}^{(j)}\geq J_{D_2}^{(j)}$.

We shall also need the continuity property of the functions just
introduced with respect to the increasing family of domains.

\proclaim{Proposition 1}

(1) Let $D$ be a bounded domain in $\Bbb C^n$. Let $D=\bigcup\sb{\nu=1}\sp{\infty}D_{\nu}$ where
$D_{\nu}\subset D_{\nu+1}$, $D_{\nu}$ is a domain in $\Bbb C^n$.
Then for any $j$ the sequence $(J_{D_{\nu}}^{(j)})_{\nu}$ is
increasing and convergent locally uniformly on $D\times\Bbb C^n$ to
$J_{D}^{(j)}$. In particular, the sequence $(\beta_{D_{\nu}})$
(respectively, $(R_{D_{\nu}})_{\nu}$) is locally uniformly
convergent to $\beta_D$ (respectively, $R_D$) on $D\times(\Bbb
C^n\setminus\{0\})$.

(2) Let $D$ be a bounded domain in $\Bbb C^n$. 
Assume that $D=\bigcup\sb{\nu=1}\sp{\infty}G_{\nu}$ where $G_{\nu}$ i a domain in $\Bbb C^n$. 
Assume additionally that for any compact set $K\subset
D$ there is a $\nu_0$ such that $K\subset G_{\nu}$ for any $\nu\geq
\nu_0$. Then the sequence $(J_{G_{\nu}}^{(j)})\sb{\nu=1}\sp{\infty}$
is locally uniformly conergent to $J_D^{(j)}$. In particular, the
sequence $(\beta_{G_{\nu}})$ (repectively, $(R_{G_{\nu}})$) is
locally uniformly convergent to $\beta_D$ (respectively, $R_D$) on
$D\times(\Bbb C^n\setminus\{0\})$.
\endproclaim

For a domain $D\subset\Bbb C$, $z\in D$ we put
$J_{D}^{(j)}(z):=J_D^{(j)}(z;1)$, $\beta_D(z):=\beta_D(z;1)$, 
$R_D(z):=R_D(z;1)$. Recall that $J_D^{(j)}=J_{D\setminus A}^{(j)}$ on $D\setminus A$ where
$A$ is a closed polar set in $D$ such that $D\setminus A$ is connected.

Denote $P(\lambda_0,r,R):=\{\lambda\in\Bbb C:r<|\lambda-\lambda_0|<R\}$, $0\leq r<R\leq\infty$, $\lambda_0\in\Bbb C$. We also put $P(r,R):=P(0,r,R)$.

\comment
$$
\align
J_D^{(0)}(z):=&\sup\{|f(z)|^2:f\in L_h^2(D),||f||_D\leq 1\}\\
J_D^{(1)}(z):=&\sup\{|f^{\prime}(z)|^2:f\in L_h^2(D),||f||_D\leq 1, f(z)=0\},\\
J_D^{(2)}(z):=&\sup\{|f^{\prime\prime}(z)|^2:f\in L_h^2(D),||f||_D\leq 1,
f(z)=0,f^{\prime}(z)=0\}
\endalign
$$
($||\cdot||_D$ denotes the $L^2$ norm on $D$).
\endcomment

We are going to prove the following result.

\proclaim{Theorem 2} Let $r\in(0,1)$, $\alpha\in(0,1)$. Then 
$$
\gather
r^{2\alpha}J_{P(r,1)}^{(0)}(r^{\alpha})\sim
\frac{1}{-\log r},\quad r^{4\alpha}J_{P(r;1)}^{(1)}(r^{\alpha})\sim
\frac{2r^{2\alpha}+2r^{2(1-\alpha)}}{1-r^2},\\
r^{6\alpha}J_{P(r;1)}^{(2)}(r^{\alpha})=\frac{A(r)}{B(r)}, \\
\text{ where }
A(r)\sim
\frac{r^2}{(1-r^2)^2}(-2^4)+\frac{r^{6(1-\alpha)}}{(1-r^2)(1-r^4)}(A)+\frac{r^{6\alpha}}{(1-r^2)(1-r^4)}(-2^5),\\
B(r)\sim \frac{2r^{2\alpha}+2r^{2(1-\alpha)}}{1-r^2}
\endgather
$$
for some $A<-100$. The symbol $\phi(r)\sim\psi(r)$ means that for any sufficiently small $\epsilon>0$  $\phi(r)-\psi(r)=\psi(r)o(r^{\epsilon})$.

In particular, 
$$
\align
\lim\sb{r\to 0^+}R_{P(r;1)}(r^{alpha})=-\infty & \text{ for $\alpha\in(1/3,2/3)$}\\
\lim\sb{r\to 0^+}R_{P(r;1)}(r^{\alpha})=2 & \text{ for
$\alpha\in(0,1/3]\cup[2/3,1)$},
\endalign
$$
\endproclaim

The above theorem gives a generalization of a result from \cite{Din~2009} 
(where the cases $\alpha=1/2$, $\alpha=0.3$ and $\alpha=0.7$ have been done). It gives an answer to a problem posed in \cite{Din~2009} on the asymptotic behavior of $R_{P(r;1)}(r^{\alpha})$ for arbitrary $\alpha\in (0,1)$. Additionally, we present in Remark 4 the precise asymptotic behavior of $R_{P(r;1)}(r^{\alpha})$ as $r\to 0^+$.

Analoguously as in \cite{Din~2009} we may make use of Theorem 2 to construct an infinitely connected planar bounded domain
with the supremum of the sectional curvature equal to $2$ and its infimum equal to $-\infty$. 
The domain constructed by us 
is a Zalcman-type domain (unlike that in \cite{Din~2009}) and the method of the proof of the above fact does not use, in contrast to \cite{Din~2009}, any sophisticated method.

Recall that the example from \cite{Din~2009} (and certainly also the one presented in Corollary 3) may be seen as the final one
in presenting examples where the supremum of the sectional curvature may be $2$ (see \cite{Chen-Lee~2009}) or its infimum may be equal to $-\infty$ (see \cite{Her~2009}) -- the example has similtanuously both properties.

\proclaim{Corollary 3} Let $\theta\in(0,1)$. Then there is a strictly
increasing sequence $(n_k)_k$ of positive integers such that $\bar
\triangle(\theta^{n_k},\theta^{2n_k})\cap\bar\triangle(\theta^{n_l},\theta^{2n_l})=\emptyset$,
$k\neq l$, $\bar\triangle(\theta^{n_k},\theta^{2n_k})\subset \frac{1}{2}\Bbb
D$ and
$$
\sup\{R_D(z):z\in D\}=2,\qquad\inf\{R_D(z):z\in D\}=-\infty,
$$
where $D=\frac{1}{2}\Bbb
D\setminus(\bigcup\sb{k=1}\sp{\infty}\bar\triangle(\theta^{n_k},\theta^{2n_k})\cup\{0\})$.
\endproclaim

\demo{Proof of Theorem 2} We start with the analysis of some more general situation.
For $0<r<R$ denote $\alpha_n^{r,R}:=||\lambda^n||^2_{P(r,R)}$, $n\in\Bbb Z$.

Note that
$$
\frac{1}{2\pi}\alpha_n^{r,R}=
\cases
\frac{R^{2(n+1)}-r^{2(n+1)}}{2(n+1)},&n\neq-1\\
\log R-\log r,&n=-1
\endcases.
$$

For $f\in L_h^2(P(r,R))$, $f(\lambda)=\sum\sb{n\in\Bbb Z}a_n\lambda^n$ the following identity
$$
||f||^2_{P(r,R)}=\sum\sb{n\in\Bbb Z}|a_n|^2\alpha_n^{r,R}
$$
holds.

Assume now that $r<1<R$.

Notice that
$$
\multline
|f(1)|^2=\left|\sum\sb{n\in\Bbb Z}a_n\right|^2=
\left|\sum\sb{n\in\Bbb Z}a_n\sqrt{\alpha_n^{r,R}}\frac{1}{\sqrt{\alpha_n^{r,R}}}\right|^2
\leq\\
\sum\sb{n\in\Bbb Z}|a_n|^2\alpha_n^{r,R}\sum\sb{n\in\Bbb Z}\frac{1}{\alpha_n^{r,R}}=||f||_{P(r,R)}^2\sum\sb{n\in\Bbb Z}\frac{1}{\alpha_n^{r,R}}.
\endmultline
$$

Therefore, $J_{P(r,R)}^{(0)}(1)\leq\sum\sb{n\in\Bbb Z}
\frac{1}{\alpha_n^{r,R}}$.

In fact, the equality above holds --
to see the equality it is sufficient to take $f\in L_h^2(P(r,R))$
with $a_n=\frac{1}{\alpha_n^{r,R}}$.

Our next aim is to give the formula for $J_{P(r,R)}^{(1)}(1)$ (which together with the previous one and general properties of the Bergman metric gives a formula for the Bergman metric of an arbitrary annulus at any point -- see Remark 4).

We prove the equality
$$
J_{P(r,R)}^{(1)}(1)=\sum\sb{n\in\Bbb Z}\frac{(n-\beta)^2}{\alpha_n^{r,R}},\tag{1}
$$
for suitably chosen $\beta\in\Bbb R$ (to be given precisely later).

Let us start with $f\in L_h^2(P_{r,R})$ of the form
$f(\lambda)=\sum\sb{n\in\Bbb Z}a_n\lambda^n$ such that $\sum\sb{n\in\Bbb Z}a_n=f(1)=0$.

For such an $f$ the following estimates hold
$$
\multline
|f^{\prime}(1)|^2=|\sum\sb{n\in\Bbb Z}na_n|^2=\left|\sum\sb{n\in\Bbb Z}(n-\beta)a_n\right|^2=
\left|\sum\sb{n\in\Bbb Z}\frac{n-\beta}{\sqrt{\alpha_n^{r,R}}}a_n\sqrt{\alpha_n^{r,R}}\right|^2\leq\\
\sum\sb{n\in\Bbb Z}\frac{(n-\beta)^2}{\alpha_n^{r,R}}\sum\sb{n\in\Bbb Z}|a_n|^2\alpha_n^{r,R}=\sum\sb{n\in\Bbb Z}\frac{(n-\beta)^2}{\alpha_n^{r,R}}||f||_{P(r,R)}^2.
\endmultline
$$
This gives the inequality '$\leq$' (with arbitrarily chosen $\beta$).
Now we take $f$ with $a_n=\frac{n-\beta}{\alpha_n^{r,R}}$, where $\beta$ is such that the equality $\sum\sb{n\in\Bbb Z}a_n=f(1)=0$ holds. If such a choice of $\beta$ could be made we would get the equality in \thetag{1}.
But this means that we need to find $\beta$ such that $\sum\sb{n\in\Bbb Z}\frac{n-\beta}{\alpha_n^{r,R}}=0$, which however is satisfied exactly if
$$
\beta=\frac{\sum\sb{n\in\Bbb Z}\frac{n}{\alpha_n^{r,R}}}{\sum\sb{n\in\Bbb Z}\frac{1}{\alpha_n^{r,R}}}.
$$
Consequently, with such a $\beta$ we get the equality
$$
\multline
J_{P(r,R)}^{(1)}(1)=\sum\sb{n\in\Bbb Z}\frac{(n-\beta)^2}{\alpha_n^{r,R}}=\sum\sb{n\in\Bbb Z}\frac{n^2-\beta n}{\alpha_n^{r,R}}+\beta\sum\sb{n\in\Bbb Z}\frac{\beta-n}{\alpha_n^{r,R}}=\\
\sum\sb{n\in\Bbb Z}\frac{n^2-n\beta}{\alpha_n{r,R}}=\frac{\phi_{r,R}(2)\phi_{r,R}(0)-\phi_{r,R}(1)^2}{\phi_{r,R}(0)},
\endmultline
$$
where $\phi_{r,R}(j):=\sum\sb{n\in\Bbb Z}\frac{n^j}{\alpha_n^{r,R}}$.

\bigskip

Let us now go on to the case of the annulus $P(r,1)$ where $0<r<1$. Our aim is to get the asymptotic
behaviour of the curvature of $P(r,1)$ at $r^{\alpha}$ (for a fixed $\alpha\in(0,1)$) as $r\to 0^+$.
First recall that
$$
J_{P(r,1)}^{(j)}(r^{\alpha})=r^{-2(j+1)\alpha}J_{P(r^{1-\alpha},r^{-\alpha})}(1).
$$

For simplicity we shall use the notation $\alpha_n=\alpha_n^{r^{1-\alpha},r^{-\alpha}}$ and
$J^{(j)}(1)=J_{P(r^{1-\alpha},r^{-\alpha})}(1)$. Then we get the following formulas
$$
\frac{\alpha_n}{2\pi}=
\cases
\frac{1-r^{2(n+1)}}{2(n+1)r^{2(n+1)\alpha}},&n\neq -1\\
-\log r,&n=-1
\endcases.
$$
From now on we forget about the constant $2\pi$.

Note that for $n\geq 0$ the following formula holds
$$
\alpha_{-n-2}=\frac{1-r^{2(n+1)}}{2(n+1)r^{2(n+1)(1-\alpha)}}.
$$

Let us define some functions (for $j=0,1,\ldots$)
$$
\multline
\phi(j):=\sum\sb{n\in\Bbb Z}\frac{n^j}{\alpha_n}=\\
\frac{(-1)^j}{-\log r}+
\sum\sb{n=0}\sp{\infty}\frac{2(n+1)}{1-r^{2(n+1)}}(n^jr^{2(n+1)\alpha}+(-1)^j(n+2)^jr^{2(n+1)(1-\alpha)})
=:\frac{(-1)^j}{-\log r}+\psi(j).
\endmultline
$$
Then we may write the formula we have just obtained in the following form:
$$
J^{(0)}(1)=\phi(0),\quad J^{(1)}(1)=\frac{\phi(2)\phi(0)-\phi(1)^2}{\phi(0)}.
$$
Note that the above formulas depend on $r$ and $\alpha$.

Our next aim is to find the formula for $J^{(2)}(1)$. We proceed similarly as above.

\comment
Therefore, the following formula holds
$$
J^{(0)}(r^{\alpha})=\frac{1}{-\log r}+\sum\sb{n=0}\sp{\infty}\frac{2(n+1)}{1-r^{2(n+1)}}\left(r^{2(n+1)\alpha}+r^{2(n+1)(1-\alpha)}\right).
$$

It is easy to see that for any $\epsilon>0$ 
$$
J^{(0}(r^{\alpha})=\frac{1}{-\log r}+\frac{2r^{2\alpha}}{1-r^2}+\frac{2r^{2(1-\alpha)}{1-r^2+o(r^{2\alpha+\epsilon})+o(r^{2(1-\alpha)+\epsilon}).
$$
Since
$$
J^{(1)}(r^{\alpha})=\frac{(1+\beta)^2}{-\log r}+\sum\sb{n=0}\sp{\infty}\frac{(n-\beta)^22(n+1)}{1-r^{2(n+1)}}\left(r^{2(n+1)\alpha}+r^{2(n+1)(1-\alpha)}\right),
$$
where
$$
\beta(r)=\beta=\frac{\frac{-1}{-\log r}+\sum\sb{n=0}\sp{\infty}\frac{2(n+1)}{1-r^{2(n+1)}}(nr^{2(n+1)\alpha}-(n+2)r^{2(n+1)(1-\alpha)})}
{\frac{1}{-\log r}+\sum\sb{n=0}\sp{\infty}\frac{2(n+1)}{1-r^{2(n+1)}}(r^{2(n+1)\alpha}+r^{2(n+1)(1-\alpha)})}.
$$

It is easy to see that $\beta(r)\sim -1$. It is easy to find that
$$
1+\beta(r)=\frac{\sum\sb{n=0}\sp{\infty}\frac{2(n+1)}{1-r^{2(n+1)}}(n+1)
(r^{2(n+1)\alpha}-r^{2(n+1)(1-\alpha)})}
{\frac{1}{-\log r}+\sum\sb{n=0}\sp{\infty}\frac{2(n+1)}{1-r^{2(n+1)}}(r^{2(n+1)\alpha}+r^{2(n+1)(1-\alpha)})}.
$$

Now it is easy to see that
$$
J^{(1)}(r^{\alpha})\sim r^{2\alpha}+r^{2(1-\alpha)}.
$$

Our next aim is to provide analoguous results on the bahavior of $J^{(2)}(r^{\alpha}$ as $r\to 0^+$.
\endcomment

Let us start with $f\in\Cal O(P(r^{1-\alpha},r^{-\alpha}))$ with $f(\lambda)=\sum\sb{n\in\Bbb Z}a_n\lambda^n$
such that $\sum\sb{n\in\Bbb Z}a_n=f(1)=0$ and $\sum\sb{n\in\Bbb Z}na_n=f^{\prime}(1)=0$.
Then
$$
\multline
|f^{\prime\prime}(1)|^2=\left|\sum\sb{n\in\Bbb Z}n(n-1)a_n\right|^2=|\sum\sb{n\in\Bbb Z}(n^2-\beta n-\gamma)a_n|^2=\\
\left|\sum\sb{n\in\Bbb Z}\frac{n^2-\beta n-\gamma}{\sqrt{\alpha_n}}a_n\sqrt{\alpha_n}\right|^2\leq
\sum\sb{n\in\Bbb Z}\frac{(n^2-\beta n-\gamma)^2}{\alpha_n}\sum\sb{n\in\Bbb Z}|a_n|^2\alpha_n.
\endmultline
$$
As before if we find $\beta,\gamma$ such that for $a_n=\frac{n^2-\beta n-\gamma}{\alpha_n}$ the equalities
$\sum\sb{n\in\Bbb Z}na_n=\sum\sb{n\in\Bbb Z}a_n=0$ hold then we shall have the equality
$$
J^{(2)}(1)=\sum\sb{n\in\Bbb Z}\frac{(n^2-\beta n-\gamma)^2}{\alpha_n}=\sum\sb{n\in\Bbb Z}\frac{n^2(n^2-\beta n-\gamma)}{\alpha_n}.
$$
The above properties are satisfied iff
for some $\beta,\gamma\in\Bbb R$ the equalities
$$
\cases
\sum\sb{n\in\Bbb Z}\frac{n^2-\beta n-\gamma}{\alpha_n}&=0\\
\sum\sb{n\in\Bbb Z}n\frac{n^2-\beta n-\gamma}{\alpha_n}&=0
\endcases
$$
hold.

The above is equivalent to the following system
$$
\cases
\beta\sum\sb{n\in\Bbb Z}\frac{n}{\alpha_n}+\gamma\sum\sb{n\in\Bbb Z}\frac{1}{\alpha_n}&=
\sum\sb{n\in\Bbb Z}\frac{n^2}{\alpha_n}\\
\beta\sum\sb{n\in\Bbb Z}\frac{n^2}{\alpha_n}+\gamma\sum\sb{n\in\Bbb Z}\frac{n}{\alpha_n}&=
\sum\sb{n\in\Bbb Z}\frac{n^3}{\alpha_n}
\endcases
$$
Since $\left(\sum\sb{n\in\Bbb Z}\frac{n}{\alpha_n}\right)^2-\sum\sb{n\in\Bbb Z}\frac{n^2}{\alpha_n}\sum\sb{n\in\Bbb Z}\frac{1}{\alpha_n}<0$, the above system of equations has one solution
$$
\align
\beta&=\frac{\sum\sb{n\in\Bbb Z}\frac{n^2}{\alpha_n}\sum\sb{n\in\Bbb Z}\frac{n}{\alpha_n}-
\sum\sb{n\in\Bbb Z}\frac{n^3}{\alpha_n}\sum\sb{n\in\Bbb Z}\frac{1}{\alpha_n}}
{\left(\sum\sb{n\in\Bbb Z}\frac{n}{\alpha_n}\right)^2-\sum\sb{n\in\Bbb Z}\frac{n^2}{\alpha_n}\sum\sb{n\in\Bbb Z}\frac{1}{\alpha_n}}=\frac{\phi(2)\phi(1)-\phi(3)\phi(0)}{\phi(1)^2-\phi(2)\phi(0)}\\
\gamma&=\frac{\sum\sb{n\in\Bbb Z}\frac{n}{\alpha_n}\sum\sb{n\in\Bbb Z}\frac{n^3}{\alpha_n}-
\left(\sum\sb{n\in\Bbb Z}\frac{n^2}{\alpha_n}\right)^2}
{\left(\sum\sb{n\in\Bbb Z}\frac{n}{\alpha_n}\right)^2-\sum\sb{n\in\Bbb Z}\frac{n^2}{\alpha_n}\sum\sb{n\in\Bbb Z}\frac{1}{\alpha_n}}=\frac{\phi(1)\phi(3)-\phi(2)^2}{\phi(1)^2-\phi(2)\phi(0)}
\endalign
$$

Therefore, we may write the formula
$$
\multline
J^{(2)}(1)=\phi(4)-\beta\phi(3)-\gamma\phi(2)=\\
\frac{\phi(4)\phi(1)^2-\phi(4)\phi(2)\phi(0)-2\phi(3)\phi(2)\phi(1)+\phi(3)^2\phi(0)+\phi(2)^3}{\phi(1)^2-\phi(2)\phi(0)}.
\endmultline
$$

So let us fix $\alpha\in(0,1)$. Then for any $\epsilon>0$ small enough
$$
\phi(0)=\frac{1}{-\log r}+\frac{2r^{2\alpha}}{1-r^2}+\frac{2r^{2(1-\alpha)}}{1-r^2}+o(r^{2\alpha+\epsilon})+o(r^{2(1-\alpha)+\epsilon}).
$$

The asymptotic behaviour of $\phi(1)^2-\phi(2)\phi(0)$ is the following. The coefficients of the term of highest order (i.e. of $\frac{1}{(-\log r)^2}$) vanish and the term at $\frac{1}{-\log r}$ is the following
$$
-(\psi(2)+\psi(0)+2\psi(1))=-\sum\sb{n=0}\sp{\infty}\frac{2(n+1)^3}{1-r^{2(n+1)}}(r^{2(n+1)\alpha}+r^{2(n+1)(1-\alpha)}).
$$

The remaining terms are the following $\psi(1)^2-\psi(2)\psi(0)$. Therefore, one may easily verify that the asymptotic behaviour is the following. For any $\epsilon>0$ small enough
$$
\phi(1)^2-\phi(2)\phi(0)=\frac{1}{-\log r}\left(\frac{2r^{2\alpha}}{1-r^2}+\frac{2r^{2(1-\alpha)}}{1-r^2}\right)+o(r^{2\alpha+\epsilon})
+o(r^{2(1-\alpha)+\epsilon}).
$$

\comment
If $\alpha<1/2$ then $\phi(1)^2-\phi(2)\phi(0)\sim -2 r^{2\alpha}$, if $\alpha>1/2$ then
$\phi(1)^2-\phi(2)\phi(0)\sim -2(2^2+1-2\cdot2^1) r^{2\alpha}$. And for $\alpha=1/2$ we have
$\phi(1)^2-\phi(2)\phi(0)\sim -2(1+1) r$.
\endcomment

We are remained with the asymptotic behavior of $\phi(4)\phi(1)^2-\phi(4)\phi(2)\phi(0)-2\phi(3)\phi(2)\phi(1)+\phi(3)^2\phi(0)+\phi(2)^3$.

First note that the coefficients of the terms $\frac{1}{(-\log r)^j}$, $j=2,3$ vanish. On the other hand the coefficient of the term
$\frac{1}{-\log r}$ is the following
$$
\multline
\psi(1)^2-2\psi(1)\psi(4)-\psi(4)\psi(2)-\psi(4)\psi(0)-\psi(2)\psi(0)\\-
2(-\psi(2)\psi(1)+\psi(3)\psi(1)-\psi(2)\psi(3))
+\psi^2(3)-2\psi(0)\psi(3)+3\psi(2)^2.
\endmultline
$$
Let us deal with the asymptotic behaviour (as $r\to 0$) of the last expression. One may calculate 
that for any $\epsilon>0$ small enough the last expression equals
$$
\multline
\frac{r^2}{(1-r^2)^2}(-2^4)+\frac{r^{6(1-\alpha)}}{(1-r^2)(1-r^4)}(A)+\\
\frac{r^{6\alpha}}{(1-r^2)(1-r^4)}(-2^5)+
o(r^2)+o(r^{6(1-\alpha)+\epsilon})+o(r^{6\alpha+\epsilon})
\endmultline
$$
where $A<-100$.

Combining all the obtained results we easily get the desired asymptotic behavior as claimed in the theorem.
\qed
\enddemo

\subheading{Remark 4} Recall the formula for the curvature
$$
R_{P(r,1)}(r^{\alpha})=2-R(r,\alpha):=2-\frac{{J^{(0)}(1)J^{(2)}(1)}}{\left(J^{(1)}(1)\right)^2}
$$

Then the result of Theorem 2 gives, in particular, the asymptotic behavior of the expression $R(r,\alpha)$ 
(and consequently the asymptotic behaviour of the holomorphic curvature) as $r\to 0^+$ which looks as follows 
$$
\cases \frac{1}{-\log r} & \text{ for $\alpha\in(0,1/3]$}\\
\frac{1}{r^{6\alpha-2}(-\log r)} & \text{ for $\alpha\in(1/3,1/2]$}\\
\frac{1}{r^{6(1-\alpha)-2}(-\log r)} & \text{ for $\alpha\in(1/2,2/3)$}\\
\frac{1}{-\log r} & \text{ for $\alpha\in[2/3,1)$}.
\endcases
$$

\subheading{Remark 5} Note that in the proof of Theorem 2 we have obtained a formula for the Bergman kernel and metric in the annulus (compare \cite{Her~1983}, \cite{Jar-Pfl~1993}) and a relatively simple expression for the sectional curvature of the annulus. 

\comment And now for $\alpha=1/2$ the coefficient of the term $r^2$
does not vanish which means that in this case the last expression
behaves like $r^2$ and consequently the curvature behaves like
$\frac{1}{-r\log r}$. Now fix $\alpha\in(1/2,1)$. The constant of
the term $r^{4(1-\alpha)}$ is zero. The constant of the term
$r^{6(1-\alpha)}$ is less than $-100$, which together with the fact
that constant at $r^2$ is larger than $-100$ gives the following two
possibilities. If $\alpha\in(1/2,2/3]$ then the last term behaves
like $r^2$ so in this case the curvature behaves like
$\frac{1}{r^{4-6\alpha}(-\log r)}$ whereas for $\alpha\in[2/3,1)$
the behavior is like $\frac{1}{-\log r}$.

The symmetry with respect to $\alpha=1/2$ finishes the proof.

Our next aim is to prove the following result.

\proclaim{Theorem} Let $\theta\in(0,1)$. Then there is a strictly
increasing sequence $(n_k)_k$ of positive integers such that $\bar
\triangle(\theta^{n_k},\theta^{2n_k})\cap\bar\triangle(\theta^{n_l},\theta^{2n_l})=\emptyset$,
$k\neq l$, $\bar\triangle(\theta^{n_k},\theta^{2n_k})\subset 1/2\Bbb
D$ and
$$
\sup\{R_D(z):z\in D\}=2,\;\inf\{R_D(z):z\in D\}=-\infty,
$$
where $D=1/2\Bbb
D\setminus(\{0\}\cup\bigcup\sb{k=1}\sp{\infty}\bar\triangle(\theta^{n_k},\theta^{2n_k}))$.
\endproclaim
Note that it follows trivially from Wiener's criterion that $D$ from
the theorem above is hyperconvex.
\endcomment

\demo{Proof of Corollary 3} We construct
inductively sequences $(n_k)$, $(x_k)$, $(y_k)$ and $(r_k)$ such
that $\theta^{n_1}+\theta^{2n_1}<x_1,y_1<1/2$ and for any
$k=1,2,\ldots$ the following properties hold:
$\theta^{n_{k+1}}+\theta^{2n_{k+1}}<x_{k+1},y_{k+1}<\theta^{n_k}-\theta^{2n_k}$,
$\theta^{n_{k+1}}+\theta^{2n_{k+1}}<r_{k+1}<\theta^{n_k}-\theta^{2n_k}$
and for any compact $L\subset \bar\triangle(0,r_{k+1})$ for which
$\Omega=\frac{1}{2}\Bbb
D\setminus\left(\bigcup\sb{j=1}\sp{k}\bar\triangle(\theta^{n_j},\theta^{2n_j})\cup
L\right)$ is connected the inequalties $R_{\Omega}(x_j)>2-1/j$,
$R_{\Omega}(y_j)<-j$ hold for any $j=1,\ldots,k+1$.

Then we put $D:=\frac{1}{2}\Bbb
D\setminus\left(\bigcup\sb{j=1}\sp{\infty}\bar\triangle(\theta^{n_j},\theta^{2n_j})
\cup\{0\}\right)$. The properties we assumed ensure us that the domain $D$ satisfies the inequalities
$R_D(x_k)>2-1/k$, $R_D(y_k)<-k$ which finishes the proof.

We go on to the construction of the above sequences. We put
$r_1:=1/4$. The possibility of the choice of $n_1$, $x_1$, $y_1$ as
desired follows from Theorem 2 together with the biholomorphic
invariance of the sectional curvature (we have to choose $n_1$
sufficiently large). The possibility of the choice of $r_2$ follows from Proposition 1. 
Now assume the system as above has been chosen
for $j=1,\ldots,k$ (with the choice of $n_j$, $x_j$, $y_j$, $j=1,\ldots,k$ and $r_j$, $j=1,\ldots,k+1$).

First note that choosing $n_{k+1}>n_k$ so that $\theta^{n_{k+1}}+\theta^{2n_{k+1}}<r_{k+1}$ we get that the recursively defined set $D_{k+1}=\frac{1}{2}\Bbb D\setminus\left(\bigcup\sb{j=1}\sp{k+1}\bar\triangle(\theta^{n_j},\theta^{2n_j})\right)$ satisfies the property $R_{D_{k+1}}(x_j)<-j$, $R_{D_{k+1}}(y_j)>2-\frac{1}{j}$, $j=1,\ldots,k$. Moreover, notice that after we choose $n_{k+1}$ and $x_{k+1},y_{k+1}$ with $\theta^{n_{k+1}}+\theta^{2n_{k+1}}<x_{k+1},y_{k+1}<\theta^{n_k}-\theta^{2n_k}$ and $R_{D_{k+1}}(x_{k+1})>2-1/(k+1)$,
$R_{D_{k+1}}(y_{k+1})<-(k+1)$ we easily get the existence of the desired $r_{k+2}$ from Proposition 1. 
Therefore, what we need is to choose $n_{k+1}>>n_k$ and properly chosen $x_{k+1}$, $y_{k+1}$.
We choose $x_{k+1}$, $y_{k+1}$ to be equal to $\theta^{n_{k+1}}+\theta^{\alpha 2n_{k+1}}$, where
$\alpha=\frac{1}{2}$ in the case of $x_{k+1}$ and $\alpha=\frac{1}{4}$ in the case of $y_{k+1}$.

We note that the following property holds:

For any small $\epsilon>0$, $\alpha\in(0,1)$, $j=0,1,2$ and for any $s\in(0,1)$ there is an $0<r_0<s$ such that for any $0<r<r_0$
$$
1\leq\frac{J_{P(r,s)}^{(j)}(r^{\alpha})}{J_{P(r,1)}^{(j)}(r^{\alpha})}\leq r^{-\epsilon}. \tag{2}
$$ 
Actually, the left inequality is trivial. The right inequality can be proven as follows.
First note that
$$
J_{P(r,s)}^{(j)}(r^{\alpha})=J_{P(\frac{r}{s},1)}^{(j)}\left(\frac{r^{\alpha}}{s}\right)s^{-2(j+1)}.
$$
Since $\frac{r^{\alpha}}{s}=\left(\frac{r}{s}\right)^{\frac{\alpha\log r-\log s}{\log r-\log s}}$, the desired property follows from Theorem 2.

Note that
$$
P(\theta^{n_{k+1}},\theta^{2n_{k+1}},r_{k+1}-\theta^{n_{k+1}})\subset \frac{1}{2}\Bbb D\setminus
\left(\bigcup\sb{j=1}\sp{k+1}\bar\triangle(\theta^{n_j},\theta^{2n_j})\right)\subset P(\theta^{n_{k+1}},\theta^{2n_{k+1}},1).
$$
\comment
And now the desired property is shown as follows
$$
J_{P(r,s)}^{(j)}(r^{\alpha})=J_{P(\frac{r}{s},1)}^{(j)}(\frac{r^{\alpha}}{s})s^{-2(j+1)}.
$$
\endcomment
Making use of \thetag{2}, Theorem 2 and the above inclusions we get the existence of $n_{k+1}$ as claimed.

\qed
\enddemo

\subheading{Remark 6} It would be interesting to find a precise description of Zalcman type domains having the property as
stated in Corollary 3. Note that such a description (complete or at least partial) has been done for a description of the boundary behavior of the Bergman kernel, Bergman metric or Bergman completeness (see \cite{Juc~2004}, \cite{Pfl-Zwo~2003}, \cite{Zwo~2002})

The construction presented in Corollary 3 is similar to the one presented in \cite{Jar-Pfl-Zwo~2000} where the first example of a fat bounded planar domain which is not Bergman exhaustive has been presented.

\enddocument

Recall that
$$
\align
\sum\sb{n\in\Bbb Z}\frac{1}{\alpha_n}&=\frac{1}{-\log r}+\sum\sb{n=0}\sp{\infty}\frac{2(n+1)}{1-r^{2(n+1)}}(r^{2(n+1)\alpha}+r^{2(n+1)(1-\alpha)})\\
\sum\sb{n\in\Bbb Z}\frac{n}{\alpha_n}&=\frac{-1}{-\log r}+\sum\sb{n=0}\sp{\infty}\frac{2(n+1)}{1-r^{2(n+1)}}(nr^{2(n+1)\alpha}-(n+2)r^{2(n+1)(1-\alpha)})\\
\sum\sb{n\in\Bbb Z}\frac{n^2}{\alpha_n}&=\frac{1}{-\log r}+\sum\sb{n=0}\sp{\infty}\frac{2(n+1)}{1-r^{2(n+1)}}(n^2r^{2(n+1)\alpha}+(n+2)^2r^{2(n+1)(1-\alpha)})\\
\sum\sb{n\in\Bbb Z}\frac{n^3}{\alpha_n}&=\frac{-1}{-\log r}+\sum\sb{n=0}\sp{\infty}\frac{2(n+1)}{1-r^{2(n+1)}}(n^3r^{2(n+1)\alpha}-(n+2)^3r^{2(n+1)(1-\alpha)})
\endalign
$$

Our next aim is to prove the following result.

\proclaim{Theorem} There is a Zalcmann-type domain $D:=\Bbb D\setminus\bigcup\sb{k=1}\sp{\infty}\bar\triangle(\frac{1}{2^{n_k}},\frac{1}{4^{n_k}})$
such that $\inf\{R_D(z):z\in D}=-\infty$ whereas $\sup\{R_D(z):z\in D\}=2$, where $(n_k)_k$
is a strictly increasing sequence of positive integers. In particular, $D$ is hyperconvex.
\endproclaim

\enddocument